\documentclass[11pt,reqno]{article}
\usepackage{amsmath,amsthm,amsfonts,amssymb,amscd}
\usepackage[latin1]{inputenc}

\usepackage{amssymb}
\usepackage{amsthm}
\usepackage{graphicx}
\usepackage{amscd}
\newtheorem{Theorem}{Theorem}[section]

\newtheorem{Proposition}[Theorem]{Proposition}
\newtheorem{Lemma}[Theorem]{Lemma}

\theoremstyle{Definition}
\newtheorem{Definition}[Theorem]{Definition}
\newtheorem{Example}[Theorem]{Example}

\theoremstyle{Remark}
\newtheorem{Remark}[Theorem]{Remark}

\def\leaderfill{\leaders\hbox to .8em{\hss .\hss}\hfill}
\def\_#1{{\lower 0.7ex\hbox{}}_{#1}}
\def\esima{${}^{\text{\b a}}$}
\def\esimo{${}^{\text{\b o}}$}
\font\bigbf=cmbx10 scaled \magstep1

\def\A{{\mathcal{A}}}
\def\sa{{\mathcal{S}}}
\def\P{{\mathcal{P}}}
\def\C{{\mathcal{C}}}
\def\L{{\mathcal{L}}}
\def\G{{\mathcal{G}}}

\def\fa{{\mathcal{F}}}
\def\O{{\mathcal{O}}}
\def\eR{{\mathcal{R}}}
\def\M{{\mathcal{M}}}
\def\D{{\mathcal{D}}}
\def\U{{\mathcal{U}}}
\def\B{{\mathcal{B}}}
\def\E{{\mathcal{E}}}
\def\H{{\mathcal{H}}}
\def\po{{\partial}}
\def\ro{{\rho}}
\def\te{{\theta}}
\def\Te{{\Theta}}
\def\om{{\omega}}
\def\Om{{\Omega}}
\def\vr{{\varphi}}
\def\ga{{\gamma}}
\def\Ga{{\Gamma}}
\def\la{{\lambda}}
\def\La{{\Lambda}}
\def\ov{\overline}
\def\al{{\alpha}}
\def\ve{{\varepsilon}}
\def\lg{{\langle}}
\def\rg{{\rangle}}
\def\lv{{\left\vert}}
\def\rv{{\right\vert}}
\def\be{{\beta}}

\def\bh{{\mathbb{H}}}
\def\bp{{\mathbb{P}}}
\def\ee{{\mathbb{E}}}
\def\re{{\mathbb{R}}}
\def\bz{{\mathbb{Z}}}
\def\bq{{\mathbb{Q}}}
\def\bd{{\mathbb{D}}}
\def\bn{{\mathbb{N}}}
\def\bk{{\mathbb{K}}}

\def\Sep{\operatorname{{Sep}}}
\def\Re{\operatorname{{Re}}}
\def\Mon{\operatorname{{Mon}}}
\def\SL{\operatorname{{SL}}}
\def\Res{\operatorname{{Res}}}
\def\Fol{\operatorname{{Fol}}}
\def\tr{\operatorname{{tr}}}
\def\dim{\operatorname{{dim}}}
\def\Aut{\operatorname{{Aut}}}
\def\GL{\operatorname{{GL}}}
\def\Aff{\operatorname{{Aff}}}
\def\Hol{\operatorname{{Hol}}}
\def\loc{\operatorname{{loc}}}
\def\Inv{\operatorname{{Inv}}}
\def\Ker{\operatorname{{Ker}}}
\def\mI{\operatorname{{Im}}}
\def\Dom{\operatorname{{Dom}}}
\def\Id{\operatorname{{Id}}}
\def\Tni{\operatorname{{Int}}}
\def\supp{\operatorname{{supp}}}
\def\Diff{\operatorname{{Diff}}}
\def\sing{\operatorname{{sing}}}
\def\Sing{\operatorname{{sing}}}
\def\codim{\operatorname{{codim}}}
\def\grad{\operatorname{{grad}}}
\def\Ind{\operatorname{{Ind}}}
\def\deg{\operatorname{{deg}}}
\def\rank{\operatorname{{rank}}}
\def\sep{\operatorname{{Sep}}}
\def\Sat{\operatorname{{Sat}}}
\def\signal{\operatorname{{signal}}}

\renewcommand{\thefootnote}

\title{Transversely projective holomorphic foliations with singularities}

\author{Bruno Sc\'ardua}

\date{}

\begin{document}

\maketitle

\begin{abstract}
In this paper we study the classification holomorphic foliations
with singularities. The main hypothesis is the existence of a
projective transverse structure outside of an analytic invariant
subset of codimension one. We prove   classification results for
germs of foliations and for  foliations in $\ov{\bf C}\times \ov{\bf
C}$ in terms of logarithmic  and  Riccati foliations. Our main
result reads as follows: {\it Let $\fa$ be a foliation on $\ov{\bf
C} \times \ov{\bf C}$ and with a projective transverse structure in
the complement of  an invariant algebraic curve $\Lambda \subset
\ov{\bf C} \times \ov{\bf C}$. Assume that the singularities of
$\fa$ in $\Lambda$ are non-resonant generalized curves. Then  $\fa$
is a logarithmic foliation or it is a rational pull-back of a
Riccati foliation.}
\end{abstract}



\footnote{The author wishes to thank Professor César Camacho for
many valuable conversations, for coauthoring various results in this
paper and for his constant and generous support.} \tableofcontents

\def\leaderfill{\leaders\hbox to .8em{\hss .\hss}\hfill}
\def\_#1{{\lower 0.7ex\hbox{}}_{#1}}
\def\esima{${}^{\text{\b a}}$}
\def\esimo{${}^{\text{\b o}}$}
\font\bigbf=cmbx10 scaled \magstep1

\def\A{{\mathcal{A}}}
\def\sa{{\mathcal{S}}}
\def\P{{\mathcal{P}}}
\def\C{{\mathcal{C}}}
\def\L{{\mathcal{L}}}
\def\G{{\mathcal{G}}}
\def\fa{{\mathcal{F}}}
\def\O{{\mathcal{O}}}
\def\eR{{\mathcal{R}}}
\def\M{{\mathcal{M}}}
\def\D{{\mathcal{D}}}
\def\U{{\mathcal{U}}}
\def\B{{\mathcal{B}}}
\def\E{{\mathcal{E}}}
\def\T{{\mathcal{T}}}

\def\po{{\partial}}
\def\ro{{\rho}}
\def\te{{\theta}}
\def\Te{{\Theta}}
\def\om{{\omega}}
\def\Om{{\Omega}}
\def\vr{{\varphi}}
\def\ga{{\gamma}}
\def\Ga{{\Gamma}}
\def\la{{\lambda}}
\def\La{{\Lambda}}
\def\ov{\overline}
\def\al{{\alpha}}
\def\ve{{\varepsilon}}
\def\lg{{\langle}}
\def\rg{{\rangle}}
\def\lv{{\left\vert}}
\def\rv{{\right\vert}}
\def\be{{\beta}}

\def\bh{{\mathbb{H}}}
\def\bp{{\mathbb{P}}}
\def\ee{{\mathbb{E}}}
\def\re{{\mathbb{R}}}
\def\bz{{\mathbb{Z}}}
\def\bq{{\mathbb{Q}}}
\def\bd{{\mathbb{D}}}

\def\bn{{\mathbb{N}}}
\def\bk{{\mathbb{K}}}

\def\Re{\operatorname{{Re}}}
\def\SL{\operatorname{{SL}}}
\def\Res{\operatorname{{Res}}}
\def\Fol{\operatorname{{Fol}}}
\def\tr{\operatorname{{tr}}}
\def\dim{\operatorname{{dim}}}
\def\Aut{\operatorname{{Aut}}}
\def\GL{\operatorname{{GL}}}
\def\Aff{\operatorname{{Aff}}}
\def\Hol{\operatorname{{Hol}}}
\def\loc{\operatorname{{loc}}}
\def\Ker{\operatorname{{Ker}}}
\def\mI{\operatorname{{Im}}}
\def\Dom{\operatorname{{Dom}}}
\def\Id{\operatorname{{Id}}}
\def\Tni{\operatorname{{Int}}}
\def\supp{\operatorname{{supp}}}
\def\Diff{\operatorname{{Diff}}}
\def\sing{\operatorname{{sing}}}
\def\Sing{\operatorname{{sing}}}
\def\codim{\operatorname{{codim}}}
\def\grad{\operatorname{{grad}}}
\def\Ind{\operatorname{{Ind}}}
\def\deg{\operatorname{{deg}}}
\def\rank{\operatorname{{rank}}}
\def\sep{\operatorname{{sep}}}
\def\Sat{\operatorname{{Sat}}}
\def\sep{\operatorname{{sep}}}
\def\signal{\operatorname{{signal}}}
\def\ord{\operatorname{{ord}}}

\section{Transversely projective foliations with singularities}
 Let $M$ be  a complex surface. A one-dimensional  holomorphic
{\it foliation with singularities} on $M$ is a pair
$\fa=(\fa_0,\sing(\fa))$ where $\sing(\fa)\subset M$ is a discrete
set and $\fa_0$ is a holomorphic foliation in the open set $M
\setminus \sing(\fa)$. It is natural to assume that there is no
extension of $\fa_0$ to a point in $\sing(\fa)$. We call
$\sing(\fa)$ the {\it singular set} of $\fa$. By definition the
leaves of $\fa$ are the leaves of $\fa_0$. By a standard application
of Hartogs' extension theorem a foliation is given in a small
neighborhood of a singularity by a holomorphic one-form with a
singularity at the given singularity.

From now on in this paper by  {\it foliation} we shall mean a
holomorphic foliation with singularities in  a complex dimension two
space. A foliation $\fa$ is called {\it transversely projective} if
the underlying ``non-singular" foliation
$\fa_0=:\fa\big|_{M\setminus \sing(\fa)}$ is transversely
projective.
 This means that there is an open cover
$\bigcup\limits_{j\in J} U_j = M\setminus \sing(\fa)$ such that in each $U_j$ the
foliation is given by a submersion $f_j\colon U_j \to \ov{\bf C}$ and if
$U_i \cap U_j \ne \emptyset$ then we have $f_i = f_{ij}\circ f_j$ in
$U_i \cap U_j$ where $f_{ij}\colon U_i \cap U_j \to \SL(2,{\bf C})$ is
locally constant. Thus,  on each intersection $U_i \cap U_j \ne
\emptyset$, we have $f_i =
\frac{a_{ij}f_j+b_{ij}}{c_{ij}f_j+d_{ij}}$ for some locally constant
functions $a_{ij}, b_{ij}, c_{ij}, d_{ij}$ with $a_{ij}d_{ij} -
b_{ij}c_{ij} = 1$. Basic references for transversely affine and
transversely projective foliations (in the nonsingular case) are
found in \cite{Godbillon}.

 As observed in \cite{Scardua1} {\it the singularities of a foliation
admitting a projective transverse structure are all of type $df=0$
for some local meromorphic function.} In this work we
will be considering foliations which are transversely projective in
the complement of {\it codimension one invariant divisors}. Such
divisors may, a priori, exhibit   singularities which do not admit
meromorphic first integrals.

\vglue.1in Next we introduce our main model.

\begin{Example}   [Riccati Foliations, cf. \cite{Scardua1} Example 1.1 page 190]
\label{Example:Riccati} {\rm We fix affine coordinates $(x,y) \in
{\bf C}^2$ and consider a polynomial vector field $X(x,y) =
p(x)\,\frac{\po}{\po x} + \big(a(x)y^2 + b(x)y +
c(x)\big)\,\frac{\po}{\po y}$ on ${\bf C}^2$. Then $X$ defines a
{\it Riccati foliation\/} on $\ov{\bf C} \times \ov {\bf C}$ as
follows: if we change coordinates via $u = \frac 1x$\,, \, $v =
\frac 1y$ then we obtain $X(x,v) = p(x)\,\frac{\po}{\po x} -
\big(a(x) + b(x)v + c(x)v^2\big) \, \frac{\po}{\po v}\,\cdot$
Similarly for
\begin{align*}
&X(u,y) = u^{-n} [\tilde p(u)\, \frac{\po}{\po u} + \big(\tilde
a(u)y^2 + \tilde b(u)y +
\tilde c(u)\big)\frac{\po}{\po y}] \,\,\text{and}\\
&X(u,v) = u^{-n} [\tilde p(u)\, \frac{\po}{\po u} + \big(\tilde
a(u) + \tilde b(u)v + \tilde c(u)v^2\big)\frac{\po}{\po v}]
\end{align*}
The similarity of these four expressions
 shows that $\Omega$ defines a holomorphic foliation $\mathcal R$ with
isolated singularities on $\ov{\bf C} \times \ov{\bf C}$ and having
a geometric behavior as follows:

\noindent (i) $\mathcal R$ is transverse to the fibers $\{a\} \times
\ov{\bf C}$ except for invariant fibers which are given in ${\bf C}^2$
by $\{p(x)=0\}$.

\noindent (ii) If $\La = \bigcup\limits_{j=1}^r \{a_j\} \times
\ov{\bf C}$ is the set of invariant fibers then $\mathcal R$ is
transversely projective in $(\ov{\bf C}\times\ov{\bf
C})\backslash\La$. Indeed, $\mathcal R|_{(\ov{\bf C}\times\ov{\bf
C})\backslash\La}$ is conjugate to the suspension of a
representation $\vr \colon \pi_1
(\ov{\bf C}\backslash\bigcup\limits_{j=1}^r \{a_j\}) \to \bp
\SL(2,{\bf C})$.

\noindent (iii) For a generic choice of the coefficients $a(x),
b(x), c(x), p(x) \in {\bf C}[x]$ the singularities of $\mathcal R$ on
$\ov{\bf C}\times\ov{\bf C}$ are hyperbolic, $\La$ is the only
algebraic invariant set and therefore for each singularity $q\in
\sing(\mathcal R) \subset \La$ there is a local separatrix of
$\mathcal R$ transverse to $\La$ passing through $q$.

\vglue.1in


\par Now we consider the canonical way of passing from
$\ov{\bf C}\times\ov{\bf C}$ to ${\bf C P}^2$ by a bi-rational map
$\sigma\colon \ov{\bf C}\times\ov{\bf C} \to {\bf C P}^2$ obtained
as a sequence of birational maps. The resulting foliation $\fa =
\sigma_*(\mathcal R) = (\sigma^{-1})^*(\mathcal R)$ induced by
$\mathcal R$ on ${\bf C P}^2$ has the following characteristics:

\noindent (i') $\fa$ is transversely projective in  ${\bf C
P^2}\backslash\La$ where $\La \subset {\bf C P}^2$ is the union of a
finite number of projective lines of the form
$\bigcup\limits_{j=1}^r \ov{\{x=a_j\}} \subset {\bf C P}^2$ in a
suitable affine chart $(x,y) \in {\bf C}^2 \subset {\bf C P}^2$.

\noindent (ii') For a generic choice of the coefficients of
$\Omega$, the singularities of $\fa$ in $\Lambda$ are hyperbolic
except for one single dicritical singularity $q_\infty \colon
(x=\infty, y=0) \in {\bf C P}^2$ which after one blow-up originates
a nonsingular foliation transverse to the projective line except for
a single tangency point.  }
\end{Example}

\section{Projective structures and differential forms}
\label{subsection:Projectivetriples}

Let $\fa$ be  a codimension one holomorphic  foliation with singular
set $\sing(\fa)$ of codimension $\ge 2$ on a complex manifold $M$.
The existence of a projective transverse structure for $\fa$ is
equivalent to the existence of suitable triples of differential
forms as follows (see \cite{Scardua1} Section 3, page 193):

\begin{Proposition} [\cite{Scardua1}, Proposition 1.1 page 190]
\label{Proposition:forms}  Assume that $\fa$ is given by an
integrable holomorphic one-form $\Om$ on $M$ and suppose that there
exists a holomorphic one-form $\eta$ on $M$ such that
$\text{\rm{\it(Proj.1)} } d\Om = \eta \wedge \Om$. Then $\fa$ is
transversely projective on $M$ if and only if there exists a
holomorphic one-form $\xi$ on $M$ such that $\text{\rm{\it(Proj.2)}
} d\eta = \Om \wedge \xi$ and $\text{\rm{\it(Proj.3)} }d\xi = \xi
\wedge \eta$.
\end{Proposition}
\smallskip
Notice that also $\xi$ defines a foliation with a projective
transverse structure given by the triple $(\xi, -\eta, \Om)$; we
will usually denote this {\it transverse foliation\/} by
$\fa^{\perp}$. \noindent This motivates the following definition:

\begin{Definition}[projective triple] {\rm Given holomorphic  one-forms (respectively,
meromorphic one-forms) $\Om$, $\eta$ and $\xi$ on $M$ we shall say
that $(\Om,\eta, \xi)$ is a {\it holomorphic projective triple\/}
(respectively, a {\it meromorphic projective triple\/}) if they
satisfy relations {\it(Proj.1)}, {\it(Proj.2)} and {\it(Proj.3)}
above. }
\end{Definition}

With this notion Proposition~\ref{Proposition:forms} says that $\fa$
is transversely projective on $M$ if and only if the holomorphic
pair ($\Om$, $\eta$) may be completed to a holomorphic projective
triple. According to
\cite{Scardua1} we may  perform modifications in a  projective triple as follows:

\begin{Proposition}
\label{Proposition:modificationforms}
\begin{itemize}
\item[\rm(i)] Given a meromorphic projective triple $(\Om, \eta,
\xi)$ and meromorphic functions $g$, $h$ on $M$ we can define a new
meromorphic projective triple as follows:

{\rm {\it(Mod.1)}}\,\, $\Om' = g\,\Om$

{\rm {\it(Mod.2)}}\,\, $\eta' = \eta + \frac{dg}{g} + h\,\Om$

{\rm {\it(Mod.3)}}\,\, $\xi' = \frac 1g\,\big(\xi - dh - h\eta -
\frac{h^2}{2}\,\Om\big)$

\item[\rm(ii)] Two holomorphic projective triples $(\Om,\eta,\xi)$
and $(\Om', \eta', \xi')$ define the same projective transverse
structure for a given foliation $\fa$ if and only if  we have {\rm
{\it(Mod.1)}, {\it(Mod.2)}} and {\rm {\it(Mod.3)}} for some
holomorphic functions $g$, $h$ with $g$ non-vanishing.

\item[{\rm(iii)}]
Let $(\Om,\eta,\xi)$ and $(\Om, \eta, \xi')$ be meromorphic
projective triples. Then $\xi' = \xi +F\,\Om$ for some meromorphic
function $F$ in $M$ with $d\,\Om = -\frac 12\, \frac{dF}{F} \wedge
\Om$.

\end{itemize}

\end{Proposition}

\smallskip

\noindent This last proposition   implies that  suitable meromorphic
projective triples also define projective transverse structures.
\noindent We can rewrite condition (iii) on $F$ as $d(\sqrt {F}
\,\Om) = 0$. This implies that if the projective triples $(\Omega,
\eta, \xi)$ and $(\Omega, \eta, \xi ^\prime)$ are not identical then
the foliation defined by $\Omega$ is transversely affine outside the
codimension one analytical invariant subset $\Lambda=\{F=0\}\cup
\{F=\infty\}$. (\cite{Scardua1}).

\begin{Definition} {\rm A meromorphic projective triple
$(\Omega ', \eta ', \xi ')$ is {\it geometric} (or also {\it
true\/}) if it can be written locally as in  {\it(Mod.1)},
{\it(Mod.2)} and {\it(Mod.3)} for some (locally defined) holomorphic
projective triple $(\Om, \eta, \xi)$ and some (locally defined)
meromorphic functions.}
\end{Definition}

\smallskip

\noindent As an immediate consequence we obtain:

\begin{Proposition}
\label{Proposition:true} A geometric projective triple $(\Om',
\eta', \xi')$ defines a transversely projective foliation $\fa$
given by $\Om'$ on $M$.
\end{Proposition}

\smallskip

\begin{Example}   [Riccati Foliations - revisited]
\label{Example:Riccatiforms} {\rm

Fix affine coordinates $(x,y) \in {\bf C}^2$ and consider a
polynomial one-form $\Om = p(x)dy - (y^2\,c(x)-yb(x)-a(x))dx$.  Then
$\Omega$ defines a {\it Riccati foliation\/} $\mathcal R$ on
$\ov{\bf C}\times \ov{\bf C}$ as seen in
Example~\ref{Example:Riccati} above.  Now we study the Lie Algebra
associated to this example.  Put $\eta = 2\frac{dy}{y} +
\frac{p^\prime+b}{p}\,dx + \frac{2a}{yp}\,dx \, \text{and} \, \xi =
\frac{-2a}{y^2p^2}\,dx.$ Then $(\Om,\eta,\xi)$ satisfies the
projective relations stated  in Proposition~\ref{Proposition:forms}.
This shows that $\fa$ is transversely projective in $\ov{\bf
C}\times \ov{\bf C}$  minus the algebraic subset $\ov{\{x \in {\bf
C} \mid p(x)=0\}\times{\bf C}} \cup \ov{{\bf C}\times\{y=0\}}$. But
since in the case $a(x)\not\equiv 0$, only the subset $\Lambda =
\{p(x)=0\}\times\ov{\bf C}$ is $\fa$ invariant it follows that the
transverse projective structure extends to $\ov{\bf C}\times \ov{\bf
C}\backslash \Lambda$. Indeed according to
Proposition~\ref{Proposition:modificationforms}  if we define $g =
\frac{-1}{p(x)y}$ then $\eta^\prime = \eta + 2g\Om =
\frac{p^\prime-b+2yc}{p}\,dx \,\text{and} \, \xi^\prime = \xi -
2dg-2g\eta - 2g^2\Om = \frac{2c}{p^2}\,dx$; define a triple
$(\Om,\eta^\prime,\xi^\prime)$ holomorphic in $(\ov{\bf C}\times
\ov{\bf C})\setminus \Lambda$ which gives a projective structure for
$\fa$ in this affine set. This projective structure coincides with
the one given in $(\ov{\bf C}\times \ov{\bf C})\setminus
(\Lambda\cup\ov{\bf C}\times\{y=0\})$ by $(\Om,\eta,\xi)$. The
one-form $\eta$ is closed if and only if $a\equiv 0$. Therefore
$\fa$ is transversely affine in $\ov{\bf C}\times \ov{\bf
C}\backslash (\Lambda\cup\ov{\bf C}\times\{\ov{y=0\}})$ if the
projective line $\{y=0\}$ is invariant.
 The forms $(\Om, \eta^\prime, \xi^\prime)$ define a rational projective
triple and  {\it the projective transverse structure of the
foliation $\fa^{\perp}$ defined by $\xi$ extends from ${\bf
C}^2\backslash\La$ to  $\ov{\bf C}\times \ov{\bf C}$}. Indeed,
$\fa_\xi$ admits a rational first integral. We will see this is a
general fact, under suitable hypothesis on the singularities of the
foliation $\fa$ on $\ov{\bf C}\times\ov{\bf C}$, admitting a
projective transverse structure in the complementary of an algebraic
one dimensional invariant subset $\La \subset \ov{\bf C}\times
\ov{\bf C}$.}

\end{Example}

As a kind of converse of the above example we have:
\begin{Proposition} [\cite{Scardua1} Theorem 4.1 page 197]
\label{Proposition:xifirstintegral} Let $\fa$ be a foliation on a
bidisc $U\subset {\bf C}^2$ or a projective surface $U$ admitting a
meromorphic projective triple $(\Omega, \eta, \xi)$ defined in $U$.
If $\xi$ admits a meromorphic first integral in $U$ then $\fa$ is a
meromorphic pull-back of a Riccati foliation.
\end{Proposition}

\begin{proof}
If we write $\xi = g\,dR$ for some meromorphic  function $g$ then we
may replace the meromorphic  triple $(\Om, \eta,\xi)$ by
$(\Om',\eta',\xi')$ where $\Om' = g\Om$, \, $\eta' = \eta +
\frac{dg}{g}$ \,and\, $\xi' = \frac 1g\,\xi = dR$.  The relations
$d\Om' = \eta' \wedge \xi'$, \,\, $d\eta' = \Om' \wedge \xi'$, \,\,
$d\xi' = \xi \wedge \eta'$ imply that $\eta' = HdR$ for some
meromorphic function $H$.  Now we define $\om := \frac{H^2}{2}\,\xi'
- H\eta' + dH = \frac 12\, H^2dR + dH$ one-form such that $d\om =
-HdH \wedge dR$. On the other hand $\eta' \wedge \om = HdR \wedge dH
= -HdH \wedge dR$. Thus $d\om = \eta' \wedge \om$. We also have
$d\eta' = dH \wedge dR = (-\frac 12\, H^2dR + dH) \wedge dR = \om
\wedge \xi'$. The meromorphic  triple $(\om,\eta',\xi')$ satisfies
the projective relations $d\om = \eta' \wedge \om$, \, $d\eta' = \om
\wedge \xi'$, \, $d\xi' = \xi' \wedge \eta'$ and therefore by
Proposition~\ref{Proposition:modificationforms} (iii)  we conclude
that $\Om' = \om + F.\xi'$ for some meromorphic function $F$ such
that $d\xi' = \xi' \wedge \frac{1}{2}\frac{dF}{F}\,\cdot$ This
implies $dF \wedge dR \equiv 0$. By the classical Stein
Factorization theorem  we may assume from the beginning that $R$ has
connected fibers and therefore $dF \wedge dR \equiv 0$ implies $F =
\vr(R)$ for some one-variable meromorphic function $\vr(z) \in {\bf
C}(z)$\footnote{In the case where  $U$ is a projective manifold all
the meromorphic objects are rational and therefore $\vr(z)$ is also
a rational function.}. We obtain therefore $\Om' = -\frac 12\, H^2dR
+ dH + \vr(R)dR = = dH - (\frac 12\, H^2 - \vr(R))dR.$ If we define
a rational map $\sigma\colon {\bf C P}^2 \dashrightarrow \ov{\bf
C}\times\ov{\bf C}$ \, by \, $\sigma(x,y) = \big(R(x,y),
H(x,y)\big)$ on ${\bf C}^2$ then clearly $\Om' = \sigma^*(dy-(\frac
12\, y^2 - \vr(x))dx)$ and therefore $\fa$ is the pull-back $\fa =
\sigma^*(\eR)$ of the Riccati foliation $\eR$ given on $\ov{\bf
C}\times\ov{\bf C}$ by the rational one-form $\Omega_\vr:=dy -(\frac
12\, y^2 - \vr(x))dx$.
\end{proof}


\section{Irreducible singularities}

  Let $\omega= a(x,y) dx+ b(x,y) dy$ be a holomorphic one-form defined
in a neighborhood $0\in U\in {\bf C}^2$. We say that $0\in {\bf
C}^2$ is a $\it singular$ point of $\omega$ if $a(0,0)=b(0,0)=0$,
and  a $\it regular$ point otherwise. We say that $0\in {\bf C}^2$
is an $\it{irreducible}$ singular point of $\omega$ if the
eigenvalues $\lambda_1, \lambda_2$ of the linear part of the
corresponding dual vector field $X= -b(x,y)\frac{\partial}{\partial
x} + a(x,y)\frac{\partial}{\partial y}$ at $0\in {\bf C}^2$ satisfy
one of the following conditions:

(1) $\lambda_1.\lambda_2\neq 0$ and $\lambda_1/\lambda_2\notin
{\bq_+}$

(2) either $\lambda_1\neq 0$ and $\lambda_2= 0$, or viceversa.

In case (1) there are two invariant curves tangent to the
eigenvectors corresponding to $\lambda_1$ and $\lambda_2$. In case
(2) there is an invariant curve tangent at $0\in{\bf C}^2$ to the
eigenspace corresponding to $\lambda_1$. These curves are called
$\it{separatrices}$ of the foliation.

Suppose that $0\in {\bf C}^2$ is either a regular point or an
irreducible singularity of a foliation $\mathcal I$. Then in
suitable local coordinates $(x,y)$ in a neighborhood $0\in U \in
{\bf C}^2$ of the origin, we have the following local normal forms
for the one-forms defining this foliation (\cite{Camacho-Sad}):
\begin{itemize}

\item [ (Reg)] $dy=0$, whenever $0\in {\bf C}^2$ is a regular point of $\mathcal I$.

and whenever $0\in{\bf C}^2$ is an irreducible singularity of
$\tilde\fa$, then either

\item[(Irr.1)] $xdy - \la ydx + \omega_2(x,y) = 0$ where $\la \in
{\bf C}\backslash\bq_+$\,, \, $\omega_2(x,y)$ is a holomorphic
one-form with a zero of order $\ge 2$ at $(0,0)$. This is called
{\it non-degenerate singularity}. Such a  singularity is {\it
resonant} if $\lambda\in \mathbb Q_-$ and {\it hyperbolic} if
$\lambda \notin \mathbb R$, or

\item[(Irr.2)]
 $y^{t+1}dx - [x(1+\la y^t) + A(x,y)]dy = 0$\, , where
$\la \in {\bf C}$, \, $t \in \bn = \{1,2,3,\dots\}$ and $A(x,y)$ is a
holomorphic function with a zero  of order $\ge t+2$ at $(0,0)$.
This is called {\it saddle-node singularity}. The {\it strong
manifold} of the saddle-node is given by  $\{y=0\}$. If the
singularity admits another separatrix then it is necessarily smooth
and transverse to the strong manifold, it can be taken as the other
coordinate axis and will be called {\it central} manifold of the
saddle-node.
\end{itemize}

 Therefore, for a suitable choice of the coordinates, we have
$\{y=0\} \subset \sep(\mathcal I,U) \subset \{xy=0\}$, where
$\sep(\mathcal I,U)$ denotes the union of separatrices of $\mathcal
I$ through $0\in{\bf C}^2$.

\textsl{}

\begin{Definition} [\cite{Touzet}, Def. II.4.1]
{\rm Let $\fa$ be a germ of a holomorphic foliation at the origin
$0\in {\bf C}^2$. We say that  $\fa$ is {\it transversely projective
of moderate growth} if it admits a meromorphic projective triple
$(\Omega, \eta, \xi)$ defined in a neighborhood of the origin. }
\end{Definition}

We recall the following fundamental result from \cite{Touzet}:
\begin{Theorem}
[Touzet, \cite{Touzet} Theorem II.4.2 and Theorem II.3.1]
\label{Theorem:Touzet} A germ of irreducible singularity $\fa$ at
the origin $0\in\bf C^2$  which is of resonant type or saddle-node
type is  projective of moderate growth if, and only if, the germ is
a pull-back of a Riccati foliation  by a meromorphic map. A
non-degenerate non-resonant singularity $xdy - \lambda y dx +
\Omega_2(x,y)=0, \, \lambda \in \bf C \setminus \mathbb Q_+$,  is
analytically linearizable if and only if the corresponding foliation
$\fa$ is transversely projective in $U\setminus \sep(\fa,U)$ for
some neighborhood $U$ of the singularity.
\end{Theorem}

\begin{Remark}
\label{Remark:Touzet} {\rm
 The proof of the first part of  the above  theorem is based on
the study and classification of the Martinet-Ramis cocycles
(\cite{Martinet-Ramis1, Martinet-Ramis2}) of the singularity
expressed in terms of some classifying holonomy map of a separatrix
of the singularity. For a resonant singularity any
 of the two separatrices has a {\it classifying holonomy} ({\it i.e.},
 the analytical conjugacy class of the singularity germ is determined by
 the analytical conjugacy class of the holonomy map of the separatrix) and for a
saddle-node it is necessary to consider the strong separatrix
holonomy map. Thus we conclude that the proof given in \cite{Touzet}
works if we only assume the existence of a meromorphic projective
triple $(\Om', \eta', \xi')$ on a neighborhood $U_0$ of $
\Lambda\setminus (0,0)$, where $\Lambda \subset \sep(\fa,U)$ is any
separatrix in the resonant case, and the strong separatrix if the
origin is a saddle-node. }
\end{Remark}

\section{Separatrices and resolution of singularities}

Suppose $\fa$ is a complex one-dimensional foliation defined on an
open neighborhood $0\in U\subset {\bf C}^2$. The {\it resolution}
process of $\fa$ at $0\in {\bf C}^2$ consists of a finite number of
successive blow-ups originating a foliation with only irreducible
singularities. This process can be described as follows. The blow-up
of $\fa$ at $0\in {\bf C}^2$ is $(U_0, \pi_0, D_0, \fa_0)$ where
$\pi_0: U_0\to U$ is the usual blow-up map. Then, $U_0$ is a complex
2-manifold, $D_0=\pi_0^{-1}(0)\subset U_0$ is an embedded projective
line called the $\it{exceptional}$ $\it{ divisor}$, and the
restriction of the map $\pi_0$ to $U_0\setminus D_0$ is a
biholomorphism from $ U_0\setminus D_0$ to $U\setminus \{0\}$.
Moreover $\fa_0$ is the analytic foliation on $U_0$ obtained by
extension to $D_0$ of $(\pi_0|_{U_0\setminus D})^*\fa$, {\it i.e.},
the pull-back foliation $\pi_)^*(\fa)$. If  $D_0$ is tangent to
$\fa_0$, i.e. $D_0$ is a leaf plus a finite number of singularities,
we say that $D_0$ is {\it non-dicritical}. Otherwise,  $D_0$ is
transverse to $\fa_0$ everywhere except at a finite number of
points, singularities or tangency points of $\fa_0$ with $D_0$. In
this last case we say that $D_0$ is {\it dicritical}.

Proceeding by induction we define the step $\underbar 0$ as the
first blow-up $(U_0, \pi_0, D_0, \fa_0)$. We assume that $(U_k,
\pi_k, D_k, \fa_k)$ has been already defined, where $\pi_k: U_k\to
U$ is a holomorphic map, such that $D_k=\pi_k^{-1}(0)$ is a divisor,
union of a finite number of embedded projective lines with normal
crossing. The crossing points of $D_k$ are called $\it corners$. The
restriction of $\pi_k$ to $ U_k\setminus D_k$ is a biholomorphism
from $U_k\setminus D_k$ to $U\setminus \{0\}$. The foliation $\fa_k$
on $U_k$ is the pull-back of $\fa$ by the map $\pi_k$.   The
Resolution theorem of Seidenberg \cite{Seidenberg} guarantees that
after a finite number of blow-up´s all corners obtained in this
process will be either irreducible singular points or regular
points. As final product we get a  complex surface $\widetilde U$
and a proper holomorphic map $\pi\colon \widetilde U \to U$, which
is a finite composition of quadratic blow-ups, such that the {\it
exceptional divisor} $D = \pi^{-1}(0)$ is a normal crossing divisor
without triple points. Also   $D$ is a finite union of projective
lines $D = \cup_{j=1}^m \mathbb P_j$\,, \, $\mathbb P_j \simeq {\bf
C} P(1)$ with negative self-intersection in $\widetilde U$. The
pull-back foliation $\widetilde\fa = \pi^*(\fa)$ is a foliation with
isolated singularities, $\sing(\widetilde\fa) \subset D$, consisting
 of irreducible singularities. Any component $\mathbb P_j \subset D$ is either $\tilde
\fa$-invariant or, in the dicritical case, everywhere transverse to $\widetilde\fa$.

Let now $\fa$ be a foliation with isolated singularities on a
complex manifold $M$.  Given an analytic invariant curve
$\Lambda\subset U$ we may perform the {\it resolution of
singularities of $\fa$ in $\Lambda$} obtaining a proper holomorphic
map $\pi \colon \tilde M \to M$ and a foliation $\tilde
\fa=\pi^*(\fa)$ such that the singularities of $\tilde \fa$ in
$\pi^{-1}(\Lambda)$ are all irreducible. Denote by $\tilde
\Lambda\subset \tilde M$ the  {\it strict transform} of $\Lambda$,
defined as   $\tilde \La := \ov{\pi^{-1}(\La \setminus (\Lambda\cap
\sing(\fa))}$. Then $\tilde \Lambda$ is $\tilde \fa$-invariant.  The
{\it large transform} of $\Lambda$ is by definition
$\pi^{-1}(\Lambda)= \tilde \Lambda \cup D$, the union of the strict
transform $\tilde \Lambda$ and the {\it exceptional divisor}
$D=\pi^{-1}(\sing(\fa))$.

\vglue.1in
Consider now an arbitrary germ of an analytic foliation $\fa$ at an
isolated singularity $0\in {\bf C}^2$. A {\it separatrix} of $\fa$
at $0\in {\bf C}^2$ is the germ at $0\in {\bf C}^2$ of an
irreducible analytic curve which is invariant by $\fa$. Let $\fa|_U$
be a representative of the foliation defined in a neighborhood $U$
of $0\in {\bf C}^2$, the separatrix is the union of a leaf of
$\fa|_U$ and the singular point $0\in {\bf C}^2$. By Newton-Puiseux
parametrization theorem, if $U$ is small enough, there is an
analytic injective map $f \colon \mathbb D \to U$ from the unit disk
$\mathbb D \subset \bf C$ onto the separatrix, mapping the
origin to $0\in {\bf C}^2$, and nonsingular outside the origin $0
\in \mathbb D$. Therefore a separatrix locally has the topology of a
punctured disk.

We shall say that the separatrix is {\it resonant} if for any loop
in the punctured disk that represents a generator of the homotopy of
the leaf, the corresponding holonomy map is a resonant
diffeomorphism. Choose a holomorphic vector field $X$ which
generates the foliation $\fa|_U$, and has an isolated singularity at
$0\in {\bf C}^2$. Then, the separatrix is resonant if the loop
$\gamma$ generating the homotopy of the leaf in the separatrix
satisfies $\exp \int_\gamma \tr(DX)$ is a root of the unity.

By the resolution of singularities we conclude that a separatrix
$\Gamma$ of $\fa$ is the projection $\Gamma=\pi(\tilde \Gamma)$ of a
curve $\tilde \Gamma$ invariant by $\tilde\fa$ and transverse to
the exceptional divisor $\pi^{-1}(0)$. We shall say that $\Gamma$ is
a {\it dicritical separatrix} if $\tilde \Gamma$ meets the
resolution divisor at a non-singular point. Equivalently,
$\Gamma=\pi(\tilde \Gamma)$ is non-dicritical if $\tilde \Gamma$ is
the separatrix of some singularity of $\tilde \fa$.  \vglue.1in


\section{Extension through singularities}

The following results are proved in \cite{Camacho-Scardua2} and
imply the existence of a globally defined projective triple in the
situation we are dealing with:

\begin{Proposition}
\label{Proposition:technical}
 Let $\mathcal I$ be a holomorphic foliation in a neighborhood $V$ of the origin
 $0 \in {\bf C}^2$ given by the holomorphic one-form
 $\Omega$  admitting a meromorphic one-form $\eta$ in $V$ with
 $d\Omega= \eta \wedge \Omega$. Suppose that $\mathcal I$ has  an
 irreducible singularity at the origin and
 is  transversely projective in $U\setminus \sep(\mathcal I,U)$ for
 some neighborhood $U\subset V$ of the origin where $\mathcal I$ has an
 expression in irreducible normal form.  Then given a meromorphic-one
 form  $\xi$  defined in
 $U \setminus \sep(\mathcal I,U)$ such that
 $(\Om,\eta,\xi)$
 is a geometric projective triple in $U \setminus \sep(\mathcal I,U)$, we have:
\begin{itemize}
\item[\rm(1)] If the origin is a non-degenerate non-resonant  singularity
 then $\xi$  extends as a meromorphic one-form to $U$.

\item[\rm(2)] If $\xi$ extends as a meromorphic one-form to  $S^* = S-\{0\}$,
for some separatrix $S\subset \sep(\mathcal I,U)$ which is not a
central manifold in case the singularity is a saddle-node, then
$\xi$ extends as a meromorphic one-form  to $U$.

\end{itemize}

\end{Proposition}

This proposition and  the Globalization theorem  in
\cite{Camacho-Scardua2} give:

\begin{Proposition}
\label{Proposition:projectiveextensible} Let $\fa$ a holomorphic
foliation defined in a neighborhood $V$ of $0 \in {\bf C}^2$
 with an isolated singularity at the origin. Suppose that $\fa$ is
 transversely projective in
 $U\setminus \sep(\fa,U)$
 for some neighborhood $U\subset V$ of the origin where $\fa$ is given
 by a holomorphic one-form $\Omega$ admitting a meromorphic one-form $\eta$
 such that
 $d\Omega= \eta \wedge \Omega$ in $U$.
Given a  meromorphic-one form  $\xi$  defined in $U \setminus
\sep(\fa,U)$ such that $(\Om,\eta,\xi)$  is a geometric projective
triple, then the one-form $\xi$ extends to $U$ provided that at any
resonant separatrix $\Gamma$ the form $\xi$ extends to a
neighborhood of an annulus $A\subset \Gamma$ around the singularity.
 \end{Proposition}

\section{Extension to codimension one divisors}
\label{section:extension}

In this section we investigate the extension of meromorphic projective
triples to a codimension one divisor, invariant or not by the foliation.

\begin{Lemma}[extension through a point]
Let $(\Omega, \eta, \xi)$ be a meromorphic projective triple on a
complex surface $M^2$, and $\Lambda \subset M$ an irreducible
analytic subset of dimension one. Suppose that the triple defines a
projective transverse structure outside $\Lambda$. If  there is a
point $q\in \Lambda$ and a neighborhood $q\in U\subset M$ to which
the projective structure extends, then this projective structure
extends to  $M$.
\end{Lemma}
\begin{proof}
We consider the local case where the foliation $\fa$ is given by a
holomorphic one-form $\Om$ in an open subset $W\subset \bf C^n$
with isolated zeros and  admitting  a meromorphic one-form $\eta$ on
$W$ satisfying $d\Om = \eta\wedge \Om$. We can assume that
$\Omega$ and $\eta$ have poles in general position with respect to
$\Lambda$.

 For $U\subset W$
small enough we can find a holomorphic submersion $y \colon U \to
\bf C$  and meromorphic functions $g, \, h$ in $U$ such that
\[
\Omega = gdy, \eta = \frac{dg}{g} + h dy, \, \xi=-\frac{1}{g}\big[dh
+ \frac{h^2}{2} dy\big] + \ell g dy
\]
where
\[
d( \sqrt {\ell} g dy ) =0.
\]
Thus, $\sqrt{\ell} g = \vr (y)$ for some meromorphic function
$\vr(z)$ and therefore $\ell= \frac{\vr ^2(y)}{g ^2}$.
Hence we have

\[
\Omega = gdy, \eta = \frac{dg}{g} + h dy, \,
\xi=-\frac{1}{g}\big[dh + \frac{h^2}{2} dy\big] + \frac{\vr ^2(y)}{
g} dy
\]

We investigate under which conditions  we can write

\[
\Omega = \tilde g d \tilde y, \eta = \frac{d\tilde g } {\tilde g} +
\tilde h d \tilde y, \, \xi = - \frac{1}{\tilde g} \big[ d \tilde h
+ \frac{\tilde h ^2}{ 2} d \tilde y\big]
\]

for some suitable meromorphic functions $\tilde g, \tilde h, \tilde y$.

Imposing the above equations we obtain

\begin{equation}
\begin{cases}
\label{equation:1}
g dy = \tilde g d \tilde y \\  \frac{dg}{g} + h dy =
 \frac{d\tilde g } {\tilde g} + \tilde h d \tilde y \\
-\frac{1}{g}\big[dh + \frac{h^2}{2} dy\big] + \frac{\vr ^2(y)}{ g} dy=
 - \frac{1}{\tilde g} \big[ d \tilde h + \frac{\tilde h ^2}{ 2} d \tilde y\big]
\end{cases}
\end{equation}

We shall refer to equations in \eqref{equation:1} as {\it main
equations}. From $g dy = \tilde g d \tilde y $ we obtain $g = r(y)
\tilde g$ for some meromorphic function $r(y)$. This implies $d
\tilde y = r(y) dy$ and then $\frac{dg}{g} + h dy =   \frac{d\tilde
g } {\tilde g} + \frac{r^\prime(y)}{r(y)} dy + h dy$ so that
replacing in the second main equation we obtain $\frac{d\tilde g }
{\tilde g} + \tilde h d \tilde y =  \frac{d\tilde g } {\tilde g} +
\frac{r^\prime(y)}{r(y)} dy + h dy$ and then
$\frac{r^\prime(y)}{r(y)} dy + h dy = \tilde h d\tilde y = \tilde h
r(y) dy$. This last equation rewrites

\begin{equation}
\frac{r^\prime(y)}{r(y)}  + h  = \tilde h r (y)
\end{equation}

and the final form

\begin{equation}
\tilde h = \frac{1}{r(y)} \big[\frac{r ^\prime(y)}{r(y)} + h\big]
\end{equation}

Let us turn our attention to the third main equation. From this we
obtain

\[
\frac{1}{g}\big[dh + \big(\frac{h^2}{2} - \vr^2(y)\big)dy\big]  =
\frac{1}{\tilde g} \big[ d \tilde h + \frac{\tilde h ^2}{ 2} d
\tilde y\big]
\]

Then

\[
\frac{\tilde g}{g}\big[dh + \big(\frac{h^2}{2} -
\vr^2(y)\big)dy\big]  = d \tilde h + \frac{\tilde h ^2}{ 2} d \tilde
y
\]

\[
\frac{1}{r(y)} \big[dh + \big(\frac{h^2}{2} - \vr^2(y)\big)dy\big]
= d \tilde h + \frac{\tilde h ^2}{ 2} d \tilde y
\]

\[
\frac{1}{r(y)} \big[dh + \big(\frac{h^2}{2} - \vr^2(y)\big)dy\big]
= d \tilde h + \frac{\tilde h ^2}{ 2} r(y)  dy
\]

\[
dh + \big(\frac{h^2}{2} - \vr^2(y)\big)dy  = r(y)\big[d
\bigg(\frac{1}{r(y)}\big(\frac{r^\prime(y)}{r(y)} + h \big)\bigg) +
\frac{1}{2 r(y)^2}\big( \frac{r^\prime(y)}{r(y)} + h ) ^2 r(y)
dy\big]
\]

\[
dh + \big(\frac{h^2}{2} - \vr^2(y)\big)dy =  r(y)\big[d
\bigg(\frac{1}{r(y)}\big(\frac{r^\prime(y)}{r(y)} + h \big)\bigg) +
\frac{1}{2}\frac{1}{r(y)}\big(\frac{r^\prime(y)}{r(y)} + h ) ^2
dy\big]
\]

\[
dh + \big(\frac{h^2}{2} - \vr^2(y)\big)dy =  \frac{1}{2}\big(\frac{
r ^\prime(y) + h}{r(y)}\big) ^2 dy  -
\frac{r^\prime(y)}{r(y)}\big(\frac{r^\prime (y)}{r(y)} + h\big) dy
+ d \big(\frac{r^\prime(y)}{r(y)} + h \big)
\]

This last equation is equivalent to

\begin{equation}
\label{equation:4}
 - \vr^2(y) = -\frac{1}{2}\big(\frac{r^\prime(y)}{r(y)} \big)^2 +
 \big(\frac{r^\prime (y)}{r(y)} \big) ^\prime
\end{equation}

Let us put
\[
s(y) :=\frac{r^\prime(y)}{r(y)}
\]

Then equation \eqref{equation:4} rewrites
\begin{equation}
s^\prime - \frac{1}{2} s^2 = - \vr^2
\end{equation}

So, the original question is reduced to find conditions under which
the equation above has a holomorphic solution. This is the case, for
instance if $\vr$ is holomorphic. Now we need to return to equation
$\frac{r^\prime(y)}{r(y)}= s(y)$ and study its solutions. It is
clear from integration that there is a holomorphic solution, which
must be given by $r(y) = e^{\int s(y)dy}$,  if and only if the given
data $s(y)$ is either holomorphic or meromorphic with a simple pole
and  integral positive residue at $y=0$.

\noindent{\bf First case}. If $s(y)$ has a simple pole at $y=0$. We may
assume for simplicity that $s(y)= a/y$ for some $ a\in {\bf C}^*$. In this case
from the differential equation
$s^\prime - s^2 /2 = - \vr^2$
we obtain $\vr=\frac{\sqrt{2a - a^2}}{y}$.
Integrating $r(y)=e^{\int s(y)dy}$ we obtain $r(y)=y^a$. Since
$r(y)=g/\tilde g$ we have that $r(y)$ is holomorphic without zeros.
In particular we cannot have $a\ne 0$, contradiction.

\noindent{\bf Second case}. If $s(y)$ has a pole of order $m+1 \geq 2$ at $y=0$.
In this case we can assume that  $s(y)=a/y^{m+1}$ for some $m \geq 1$ and
integration gives $r(y)=e^{-\frac{a}{m y^m}}$ which is not meromorphic at the
origin, contradiction.

\noindent{\bf Third case}. If $s(y)$ is holomorphic at $y=0$. In
this case we write $s(y)=ay ^m$ for some $m \geq 0$. We obtain
$r(y)= e^{\frac{a}{m+1} y^{m+1}}$ which is holomorphic and
non-vanishing.

Let us now finish the proof. Because the projective structure
extends to $U$ the equation \eqref{equation:1} has a holomorphic
solution and this implies that $\vr(y)$ is holomorphic according to
the above considerations. As a consequence the one-form $\xi$ is
also holomorphic in $U$ and therefore admits a holomorphic extension
to $\Lambda\setminus[(\Omega)_\infty \cup (\eta)_\infty]$. Hence,
the projective structure extends  to
$\Lambda\setminus[(\Omega)_\infty \cup (\eta)_\infty]$ and then to
$\Lambda$.

\end{proof}

\begin{Lemma}
Let $(\Omega, \eta, \xi)$ be a meromorphic projective triple in a
complex surface $M$. Assume that the triple defines  a projective
transverse structure for $\fa$ in $M\setminus \Lambda$ for some
invariant codimension one analytic subset $\Lambda\subset M$. Let
$\xi^\prime$ be a meromorphic one-form in $M$ such that $(\Omega,
\eta, \xi^\prime)$ is also a projective triple. Then $\Lambda$ is
$\xi$-invariant if and only if it is $\xi^\prime$-invariant.
\end{Lemma}

\begin{proof}
We fix a local coordinate system $(x,y)\in U$ centered at a point
$p\in M$ such that $\fa$ is given in these coordinates by $\Omega=g
dy$ and $\Lambda$ by $\{y=0\}$. We may write $\xi^\prime= \xi + \ell
\Omega$ where $d(\sqrt{\ell} \Omega)=0$. Then we have $\ell=
\frac{\vr^ 2(y)}{g}$ for some meromorphic function $\vr(z)$. Assume
by contradiction that $\Lambda$ is not $\xi$-invariant  but
$\Lambda$ is $\xi^\prime$-invariant. We may assume that the polar
set of $\xi$ has no irreducible component contained in $\Lambda$ and
therefore $\vr(y)$ and $g$ have no poles on $\{y=0\}$. Write
$\xi^\prime= Adx + B dy$ with holomorphic coefficients $A(x,y), \,
B(x,y)$. Since $\Lambda$ is $\xi^\prime$-invariant we have $A(x,y)=y
\, A_1(x,y)$ for some holomorphic function $A_1(x,y)$. Then from
$\xi^\prime= \xi + \ell \Omega$ we get $\xi= y A_1(x,y) dx + ( B
(x,y) - \frac{\vr^2(y)}{g})dy$. Since $A_1$ and $ B (x,y) -
\frac{\vr^2(y)}{g}$ have no poles in $\{y=0\}$ we conclude from the
above expression that $\Lambda$ is $\xi$-invariant, contradiction.
\end{proof}

\begin{Lemma}[non-invariant divisor, \cite{Camacho-Scardua2}]
\label{Lemma:noninvariantextension} Let be given a holomorphic
foliation $\fa$ on a complex surface  $M$. Suppose that $\fa$ is
given by a meromorphic integrable one-form $\Omega$ which admits a
meromorphic one-form $\eta$ on $M$ such  that $d \Omega= \eta \wedge
\Omega$. If $\fa$ is transversely projective in $M\setminus \Lambda$
for  some {\it non-invariant}  irreducible analytic subset $\Lambda
\subset M$ of codimension one then $\fa$ is transversely projective
in $M$. Indeed, the projective transverse structure for $\fa$  in
$M\setminus \Lambda$ extends to $M$ as a projective transverse
structure for $\fa$.
\end{Lemma}

\begin{proof} Our argumentation is local, {\it i.e.},
we consider a small neighborhood $U$ of a  generic point $q\in
\Lambda$ where $\fa$ is transverse to $\Lambda$. Thus, since
$\Lambda$ is not invariant by $\fa$, performing changes as
$\Omega^\prime = g_1 \Omega$ and $\eta ^\prime = \eta +
\frac{dg_1}{g_1}$  we can assume that $\Omega$ and $\eta$ have poles
in general position with respect to $\Lambda$ in $U$. The existence
of a projective transverse structure for $\fa$ in $M\setminus\Lambda$ then
gives a meromorphic one-form $\xi$ in $M \setminus \Lambda$ such
$(\Omega, \eta, \xi)$ is a geometric projective triple in $M
\setminus \Lambda$. For $U$ small enough we can assume that for
suitable local coordinates $(x,y)\in U$ we have
$\Lambda \cap U = \{x=0\}$ and also
\[
\Omega = gdy, \eta = \frac{dg}{g} + h dy
\]
for some holomorphic function $g, h \colon U \to {\bf C}$ with $1/g$
also holomorphic in $U$. Then we have
\[
\xi=-\frac{1}{g}\big[dh + \frac{h^2}{2} dy\big]
\]
where
\[
d( \sqrt {\ell} g dy ) =0
\]
Thus, $\sqrt{\ell} g = \vr (y)$ for some meromorphic function
$\vr(y)$ defined for $x \ne 0$ and therefore for $x=0$. This
shows that $\xi$ extends  to $U$ as a {\it holomorphic one-form} and
then  the projective structure extends to $U$. This shows that the
transverse structure extends to $\Lambda$.
\end{proof}

\section{Germs of foliations and foliations on projective spaces}

 Let $\fa$
be a holomorphic foliation of codimension one on ${\bf C P}^2$
having singular set $\sing(\fa) \subsetneqq {\bf C P}^2$. As it is
well-known we can assume that $\sing(\fa)$ is of codimension $\ge 2$
and $\fa$ is given in any affine space ${\bf C}^2 \subset {\bf C
P}^2$ with coordinates $(x,y)$, by a polynomial one-form $\Om(x,y) =
A(x,y)dx + B(x,y)dy$ with $\sing(\fa) \cap {\bf C}^2 = \sing(\Om)$.
In particular $\sing(\fa) \subset {\bf C P}^2$ is a nonempty finite
set of points. Given any algebraic subset $\Lambda\subset {\bf C
P}^2$ of dimension one we can therefore always obtain a meromorphic
(rational) one-form $\Om$ on ${\bf C P}^2$ such that $\Om$ defines
$\fa$, \, $(\Om)_\infty$ is non-invariant and in general position
(indeed, we can assume that $(\Om)_\infty$ is any projective line in
${\bf C P}^2$). Also if we take $\eta_0 = \frac{B_x}{B}\,dx +
\frac{A_y}{A}\,dy$ then we obtain a rational one-form such that
$d\Om = \eta_0 \wedge \Om$ and with polar set given by
$(\eta_0)_\infty = \{(x,y) \in {\bf C}^2 : A(x,y) = 0\} \cup \{(x,y)
\in {\bf C}^2 : B(x,y) = 0\} \cup (\Om)_\infty\,.$ In particular,
$(\eta_0)_\infty \cap {\bf C}^2$ has order one and the ``residue" of
$\eta_0$ along any component $T$ of $(\Om)_\infty$ equals $-k$ where
$k$ is the order of $T$ as a set of poles of $\Om$. Any rational
one-form $\eta$ such that $d\Om = \eta \wedge \Om$ writes $\eta =
\eta_0 + h\Om$ for some rational function $h$. We obtain in this way
one-forms $\eta$ with appropriately  located set of poles, with
respect to $\fa$,  and applying  Proposition~\ref{Proposition:forms}
and ~\ref{Proposition:modificationforms} we obtain:

\begin{Proposition}
\label{Proposition:projectivetripleoff}  Let $\fa$ be a holomorphic
foliation on ${\bf C P}^2$. Assume that $\fa$ is transversely
projective in ${\bf C P}^2\backslash \Lambda$ for some algebraic
subset $\Lambda$ of dimension one. Then $\fa$ has a projective
triple $(\Om, \eta, \xi)$ on ${\bf C P}^2 \setminus \Lambda$  where
$\Om$ and $\eta$ are rational one-forms and $\xi$ is meromorphic on
${\bf C P}^2 \backslash \Lambda$. In particular $\xi$ defines a
transverse foliation $\fa^{\perp}$ to $\fa$ on ${\bf C P}^2
\backslash \Lambda$ having a projective transverse structure.
\end{Proposition}

This proposition admits a natural local version, {\it i.e.},  a
version for germs of foliations at the origin $0\in \bf C ^2$ where
the curve $\Lambda$  is replaced by a finite set of local branches
of separatrices of the foliation through the singularity.

We recall that a germ of a foliation singularity at the origin $0\in
{\bf C}^2$ is a {\it generalized curve} if it is non-dicritical and
exhibits no saddle-node in its resolution by blow-ups
(\cite{C-LN-S}). The generalized curve is {\it resonant} if {\it
all} singularities are of resonant type, otherwise it is called {\it
non-resonant}.

For this type of singularity we have the following extension lemma:

\begin{Lemma}
\label{Lemma:projectiveextensible} Let $\mathcal F$ be a germ of a
non-resonant generalized curve at the origin $0\in {\bf C}^2$.
Suppose
 that $ \fa$ is transversely projective in
 $U\setminus \sep(\mathcal \fa,U)$ and let $(\Omega, \eta, \xi)$ be a
 meromorphic triple in $U\setminus \sep(\fa, U)$ with $\Omega$ holomorphic
 in $U$, $\eta$ meromorphic in $U$ and
 $\xi$ meromorphic in $U\setminus \sep(\fa,U)$.
Then the one-form $\xi$ extends to $U$  as a meromorphic one-form.
 \end{Lemma}
\begin{proof} This lemma follows from
Proposition~\ref{Proposition:projectiveextensible}.
\end{proof}

We extend the notion of generalized curve in a natural way by
allowing dicritical components, but no saddle-nodes, in the
resolution process. Such singularities will be called {\it extended
generalized curves}. Let now $\fa$ be a foliation in $U$ where $U$
is either a bidisc centered at the origin, or a projective surface.
Consider an invariant subset $\Lambda\subset U$ analytic of
dimension one. Denote by $\pi \colon \tilde U \to U$ the resolution
morphism of the singularities of $\fa$ in $\Lambda$. We  say that
{\it the singularities in $\Lambda$ are {\it non-resonant} extended
generalized curves if each connected component of invariant of the
resolution divisor $\pi^{-1}(\Lambda)$ contains some non-resonant
(non-degenerate) singularity}. In a natural extension of the
arguments in the proof of Lemma~\ref{Lemma:projectiveextensible} we
obtain:

\begin{Lemma}
\label{Lemma:projectiveextensiblegeneralized} Let $\mathcal F$ be a
germ of a foliation at the origin $0\in {\bf C}^2$. Suppose  that $
\fa$ is transversely projective in
 $U\setminus \Lambda$ where $\Lambda\subset \sep(\mathcal \fa,U)$ is
 a finite set of local branches.
 Assume that  the singularity $0\in\Lambda$ is a   non-resonant
 extended generalized curves
 Let $(\Omega, \eta, \xi)$ be a
 meromorphic triple in $U\setminus \sep(\fa, U)$ with $\Omega$ holomorphic
 in $U$, $\eta$ meromorphic in $U$ and
 $\xi$ meromorphic in $U\setminus \sep(\fa,U)$.
Then the one-form $\xi$ extends to $U$  as a meromorphic one-form.
 \end{Lemma}

This lemma and Proposition~\ref{Proposition:projectivetripleoff} promptly give:

\begin{Theorem}
\label{Theorem:projectiveplane}  Let $\fa$ be a holomorphic
foliation on $U$ where $U$ is either the projective plane  $\bf C
P^2$ or a bidisc centered at the origin $0\in \bf C^2$. Assume that
$\fa$ is transversely projective in $U\setminus \Lambda$ where
$\Lambda \subset U$ is either an algebraic invariant curve in the
projective plane or a finite union of local branches of separatrices
of $\fa$ through the origin. Suppose that the singularities of $\fa$
in $\Lambda$ are non-resonant extended generalized curves.  Then
$\fa$ admits a meromorphic projective triple $(\Omega, \eta, \xi)$
defined in $U$, which defines the  projective transverse structure
in $U\setminus \Lambda$.
\end{Theorem}

\section{Monodromy}
In this section we follow original ideas from \cite{Nicolau-Paul}
in the same vein as  in \cite{Scardua2}. Let $\fa$ be a holomorphic
foliation with singularities on a complex surface $M$ and $X \subset
M$ an invariant codimension one analytic subset such that $\fa$ is
transversely projective in $M \backslash X$. According to
\cite{Godbillon} the foliation $\fa\big|_{M\backslash X}$ admits a
{\it development\/}, i.e., there is a Galoisian covering $p\colon P
\to M\backslash X$ where $p$ is holomorphic, a homomorphism $h\colon
\pi_1(M\backslash X) \to \SL(2,{\bf C})$ and a holomorphic
submersion $\Phi\colon P \to {\bf C  P}^1$ such that:
\begin{itemize}
\item[(i)] $\Phi$ is $h$-equivariant. \item[(ii)]
$p^*\big(\fa\big|_{M\backslash X}\big)$ is the foliation defined
by the submersion $\Phi$.
\end{itemize}

\medskip

\begin{Remark}{\rm The construction of the development in \cite{Godbillon} requires
the foliation to be nonsingular. In our case,  it  is not necessary to require
that $\sing(\fa) \cap (M\backslash X) = \emptyset$. Indeed, by
Hartogs' Extension Theorem any holomorphic map from $(M\backslash
X)\backslash (\sing(\fa) \cap (M\backslash X))$ to $\bf CP^1$  extends uniquely to
a holomorphic map from  $M\backslash X$ to $\bf CP^1$. Also, since $\codim _{\bf C} \sing
\fa = 2$ the inclusion $\pi_1((M\backslash X)\backslash ((M\backslash
X) \cap \sing(\fa)) \to \pi_1(M\backslash X)$ is an isomorphism.
Nevertheless, for our purposes it is enough to consider the case
$\sing(\fa) \subset X$.}
\end{Remark}

\smallskip

\noindent Using the notion of development we can introduce the
notion of {\it monodromy\/} of the projective transverse structure
of $\fa\big|_{M\backslash X}$ as follows:

Fix a base point $m_0 \in M\backslash X$ and  a  local
determination $f_{m_0}$ of the submersion $\Phi$ in a small ball
$B_{m_0}$ centered at $m_0$ (we have the following commutative
diagram)
$$
\begin{matrix}
&P& \supset &p^{-1}(B_{m_0})& \Phi\big|_{p^{-1}(B_{m_0})}\\
&p\downarrow\,\,\, &  &p\big|_{p^{-1}(B_{m_0})}\downarrow\qquad\searrow\,\,&\\
&M\backslash X& \supset
&B_{m_0}&\!\!\!\!\overset{f_{m_0}}{\longrightarrow}\quad {\bf C} P(1)
\end{matrix}
$$

\noindent Notice that $p^{-1}(B_{m_0}) =
\underset{\al\in\A}{\bigcup\!\!\!\!\!\cdot}\, U_\al$, \,\,
$p\big|_{U_\al}\colon U_\al \to B_{m_0}$ is a biholomorphism for
each $\al \in \A$.


\noindent By construction, the total space of the covering
$p\colon P \to M \backslash X$ is obtained by analytic
continuation of $f_{m_0}$ along all the elements in
$\pi_1(M\backslash X, m_0)$.

\noindent The fiber $p^{-1}(m_0)$ is the set of all local
determinations $f_{m_0}$ at $m_0$\,. We can, by the general theory
of transitive covering spaces, identify the group $\Aut(P,p)$ of
deck transformations of $p\colon P \to M\backslash X$ to the
quotient $\pi_1(M\backslash X;m_0)\big/ p_\#\pi_1(P;f_{m_0})$.
This is the {\it monodromy group\/} of $\fa\big|_{M\backslash X}$
which will be  denoted by $\text{Mon}(\fa,X)$.

\noindent The {\it monodromy map\/} is the natural projection
$$
\begin{matrix}
\ro\colon \pi_1(M\backslash X;m_0) \longrightarrow
&{\pi_1(M\backslash X;m_0)}&\!\!\!\big/p_\#\pi_1(P;f_{m_0}) =
\text{Mon}(\fa,X)&
\end{matrix}
$$
Our first remark is the following:

\medskip

\begin{Lemma} The monodromy group $\Mon(\fa,X)$ is
naturally isomorphic to a subgroup of $\SL(2,{\bf C})$.
\end{Lemma}

\smallskip

\begin{proof} This is clear since $\fa\big|_{M\backslash X}$ is
transversely projective on $M\backslash X$.
\end{proof}

\section{Holonomy versus monodromy}
\label{subsection:monodromyholonomy}

Here we keep on following arguments originally in
\cite{Nicolau-Paul}  and mimed in \cite{Scardua2}. We proceed to
study the holonomy of each irreducible component of $X$. It is
enough to assume that $X$ is the union of a smooth compact curve
$\La$ and local analytic separatrices $\sep(\fa,\La)$ of $\fa$
transverse to $\La$;\,\, $X = \La \cup \sep(\fa,\La)$, all of them
smooth invariant and without triple points. We suppose that
$\sing(\fa) \cap \Lambda \ne \emptyset$, each singular point in
$\Lambda$ is irreducible and, if it admits two separatrices then one
is transverse to $\Lambda$). In this case we can consider a
$C^{\infty}$ retraction $r\colon W \to \La$ from some tubular
neighborhood $W$ of $\La$ on $M$ onto $\La$ such that, $\forall\,m
\in \La$ the fiber $r^{-1}(m)$ is either a disc transverse to $\fa$
or a local branch of $\sep(\fa,\La)$ at $m \in \sing(\fa)$. We set
$V = W \backslash(X \cap W)$ to obtain a $C^\infty$ fibration
$r\big|_V \colon V \to \La\backslash\sing(\fa)$ by punctured discs
over $\La\backslash\sing(\fa)$. Since
$\pi_2(\La\backslash\sing(\fa)) = 0$ the homotopy exact sequence of
the above fibration gives the exact sequence
$$
0 \longrightarrow \bz \longrightarrow \pi_1(V,\tilde m_0)
\overset{\tau}{\longrightarrow} \pi_1(\La\backslash\sing(\fa);m_0)
\longrightarrow 0
$$
where $\tilde m_0 \in V$ is a base point and $m_0 \in
\La\backslash\sing(\fa)$ is its projection and
$\tau=(r\big|_V)_\#$.

\noindent Now we consider the restriction of the covering space
$P$ to $V$; indeed for our purposes we may assume that $W = M$ and
$V = M\backslash X$ so that we are just considering the space $P$
itself. Let $\ro$ be the monodromy map
$$
\begin{matrix}
\ro\colon \pi_1(V;\tilde m_0) \longrightarrow
 &{\pi_1(V;\tilde m_0)} &\!\!\!\!\!\big/p_\#(\pi_1(p^{-1}(V);
 f_{\tilde m_0}))=:\text{Mon}(\fa,V)&
\end{matrix}
$$

\noindent Denote by $\Mon(\fa,\La)$ the quotient of $\Mon(\fa,V)$
by the (normal) subgroup \linebreak $\Ker(\tau) \cong \bz$. Then
there is a unique morphism $[\ro]$ such that the diagram commutes:
$$
\begin{matrix}
0 \longrightarrow \bz &\to& \pi_1(V;\tilde m_0)& \longrightarrow
&\pi_1(\La\backslash\sing(\fa);m_0) \to 0\\
&\searrow & \,\,\,\ro\,\downarrow& &\,\,[\ro]\,\downarrow\\
& &\Mon(\fa,V)& \longrightarrow &\Mon(\fa,\La) \to 0
\end{matrix}
$$

\noindent The morphism $[\ro]$ is a monodromy of $\fa\big|_V$ seen
as follows: \newline  given any element $[\ga] \in
\pi_1(\La\backslash\sing(\fa);m_0)$  the monodromy $[\ro]([\ga])$ is
the analytic continuation of the local first integral $f_{m_0}$
along $\ga$ and its holonomy lifting. This gives:

\medskip

\begin{Lemma} There exists a surjective group
homomorphism $\al\colon \Hol(\fa,\La) \longrightarrow
\Mon(\fa,\La)$ such that the diagram commutes
$$
\begin{matrix}
&\pi_1(\La\backslash\sing(\fa))&\\
\Hol\,\,\quad\swarrow& & \searrow\,\, [\ro]\quad\\
\Hol(\fa.\La) &\overset{\al}{\longrightarrow}& \Mon(\fa;\La)
\end{matrix}
$$

\noindent where $\Hol\colon \pi_1(\La\backslash\sing(\fa))
\longrightarrow \Hol(\fa;\La)$ is the holonomy morphism of the
leaf $\La\backslash\sing(\fa)$ of $\fa$, and $[\ro]\colon
\pi_1(\La\backslash\sing(\fa)) \longrightarrow \Mon(\fa;\La)$ is
as above.
\end{Lemma}

\smallskip

\noindent The kernel of $\al$ is the subgroup $\Ker(\al) <
\Hol(\fa;\La)$ of those diffeomorphisms keeping fixed any element
$\ell(z)$ of the fiber of $p\big|_V\colon V \to
\La\backslash\sing(\fa)$ over $m_o \in \La\backslash\sing(\fa)$.
Therefore $\Ker(\al)$ is a subgroup of the {\it invariance
group\/} of $\ell$, \, $\Inv(\ell,z)$, defined as follows
$\Inv(\ell,z) = \big\{h \in \Diff({\bf C},0); \ell \circ h \equiv
\ell\big\}$, in the sense that if $p_\ell\colon V_\ell \to \bd^*$
is the covering space of the punctured disc $\bd^* =
\bd\backslash\{0\}$ associated to $\ell$ then $\ell \circ h \equiv
\ell$ {\it means that\/} $\forall\, m \in \bd^*$, $\forall\,\ell_m
\in p_\ell^{-1}(m)$, $\exists\,\ell_{h(m)} \in p_\ell^{-1}(h(m))$,
\, $\ell_{h(m)} \circ h = \ell_m$\,.

\noindent In particular, to any element $h \in \Inv(\ell,z)$ there
is associated a pair $(\tilde h,h)$ where $\tilde h$ is the
lifting of $h$ to the covering space $V_\ell$ defined by $\tilde
h\colon \ell_m \mapsto \ell_{h(m)}$\,.

\smallskip

\noindent Another lemma we need is:

\medskip

\begin{Lemma}
\label{Lemma:solvablegroups} Let $0 \to G \to H \to K \to 0$ be an
exact sequence of groups. Then $H$ is solvable if, and only if, $G$
and $K$ are solvable.
\end{Lemma}

\smallskip

\noindent From the above discussion we have an exact sequence
$$
0 \longrightarrow \Ker(\al) \longrightarrow \Hol(\fa,\La)
\overset{\al}{\longrightarrow} \Mon(\fa,\La) \longrightarrow 0
$$
We claim  that $\Inv(\ell,z)$ is solvable. Indeed, suppose the
contrary. By Nakai's Density Lemma \cite{Nakai} the orbits of a
non-solvable subgroup of $\Diff({\bf C},0)$ are locally dense in a
neighborhood $\Gamma$ of the origin. Let therefore $m\in \Gamma$ be
a point and $\Gamma_m \subset \Gamma \setminus \{0\}$ be a small
sector with vertex at the origin, such that the orbit of $m$ in
$\Gamma_m$  is dense in $\Gamma_m$. Denote by $\ell_{\Gamma_m}$ a
local determination of $\ell$ in $\Gamma_m$. Then $\ell_{\Gamma_m}$
is constant along the orbits of $\Inv(\ell,z)$ in $\Gamma_m$ and the
orbit of $m$ is dense in $\Gamma_m$ so that $\ell_{\Gamma_m}$ is
constant in $\Gamma_m$. By analytic continuation $\ell$ and the
first integral $\Phi$ are constant yielding a contradiction. Thus
the group $\Inv(\ell,z)$ is solvable and therefore embeds in
$\SL(2,{\bf C})$. Hence $\Hol(\fa,\La)\big/\Ker(\al) \simeq
\Mon(\fa,\La)$ embeds in $\SL(2,{\bf C})$ but $\Hol(\fa,\La)$ embeds
in $\Diff({\bf C},0)$, as well as $\Ker(\al)$ embeds in $\Inv(\ell)$
which is a subgroup of $\Diff({\bf C},0)$ and therefore
$\Hol(\fa,\La)\big/\Ker(\al)$ is isomorphic to a subgroup of
$\SL(2,{\bf C})$ with a fixed point. This implies that indeed,
$\Hol(\fa,\La)\big/\Ker(\al)$ is solvable and conjugate to a
subgroup of $\Aff({\bf C},0)$. Therefore $\Mon(\fa,\La)$ is solvable
and by Lemma~\ref{Lemma:solvablegroups} the holonomy group
$\Hol(\fa,\La)$ is solvable.

\par Summarizing the above discussion we have:

\begin{Theorem}
\label{Theorem:solvable}  Let $\fa$ be a holomorphic foliation on
a complex surface $M$, \, $X \subset M$ a closed analytic
invariant curve and assume that  $\fa$ is transversely projective
in $M\backslash X$. Let $\Lambda\subset X$ be an irreducible
component of $X$. We suppose that each singular point in $\Lambda$
is irreducible and exhibits a separatrix transverse to $\Lambda$.
Then the holonomy group $\Hol(\fa,\La)$ of the leaf
$\La\backslash(\sing(\fa) \cap \La)$ of $\fa$ is a solvable
group.\end{Theorem}
\smallskip
\medskip


\section{Classification of transversely projective foliations}
\label{section:applications}

We consider now an application of the above study to the
classification of foliations with projective transverse structure.

\begin{Theorem}
\label{Theorem:germs}  Let $\fa$ be a germ of holomorphic foliation
at the origin $0 \in {\bf C}^2$. Suppose that
\begin{itemize}
\item[\rm(i)] $\fa$ is a germ of a non-resonant generalized curve and can be
reduced with a single blow-up.

\item[\rm(ii)] $\fa$ is transversely projective
outside of the set $\sep(\fa,0)$ of local separatrices of $\fa$
through $0$.
\end{itemize}

\noindent Then $\fa$ is given by a logarithmic one-form in a
neighborhood of the origin or it is a meromorphic pull-back of a
germ of a Bernoulli type foliation $\eR\colon \al(x)dy -
(y^2\be_0(x) + y\be_1(x))dx = 0$ where $\al$, $\be_0$, $\be_1$,
$\be_2$ are meromorphic in some neighborhood of the origin.
\end{Theorem}

We shall need the following well-known technical result.
\begin{Lemma}
\label{Lemma:linearization} Let $G < \Diff({\bf C},0)$ be a solvable non-abelian
subgroup of germs of holomorphic diffeomorphisms fixing the origin
$0 \in {\bf C}$.
\begin{itemize}
\item[\rm(i)] If  the group of commutators
$[G,G]$ is not cyclic then $G$ is analytically conjugate to a
subgroup of $\bh_k = \big\{z \mapsto
\frac{az}{\sqrt[k]{1+bz^k}}\big\}$ for some $k \in \bn$.
\item[\rm(ii)] If there is some $f \in G$  of the form $f(z) = e^{2\pi
i\la}\,z +\dots$ with $\la \in{\bf C}\backslash \bq$ then $f$ is
analytically linearizable in a coordinate that also embeds $G$ in
$\mathbb H_k$.
\end{itemize}

\end{Lemma}

\begin{proof} (i) is in \cite{Cerveau-Moussu}. Given
$f \in G$ as in (ii) then by (i) we can write $f(z) = \frac{e^{2\pi
i\la}\,z}{\sqrt[k]{1+bz^k}}$ for some $k \in \bn$, \, $b \in {\bf C}$.
Since $\la \in {\bf C}\backslash\bq$ the homography $H(z) =
\frac{e^{2\pi i\la}\,z}{1+bz}$ is conjugate by another homography to
its linear part $z \mapsto e^{2\pi i\la}\,z$ and therefore $f$ is
analytically linearizable.
\end{proof}

\begin{proof}[Proof of Theorem~\ref{Theorem:germs}] Let
$\fa$ be given in an open subset $0 \in U \subset {\bf C}^2$ and put
$\widetilde\fa = \pi^*(\fa)$ in $\widetilde U = \pi^{-1}(U)$ where
$\pi\colon \widetilde{{\bf C}}_0^2 \to {\bf C}^2$ is the blow-up of
${\bf C}^2$ at $0 \in {\bf C}^2$. Then the exceptional divisor $\La
= \pi^{-1}(0)$ is a compact invariant curve and we have
$\sep(\widetilde\fa,\La) = \ov{\pi^{-1}(\Sep(\fa,0)\backslash\{0\})}
= \ov{\pi^{-1}(\Sep(\fa,0))\backslash\La}$ in $\widetilde U$.
Therefore, each singularity of $\tilde \fa$ in $\Lambda$ is
irreducible and exhibits a separatrix transverse to $\Lambda$. Now,
by hypothesis (ii) the pull-back foliation $\widetilde\fa$ is
transversely projective in $\widetilde U \backslash\widetilde X$
where $\widetilde X = \La \cup \Sep(\widetilde\fa,\La)$. According
to Theorem~\ref{Theorem:solvable}  this implies that the holonomy
group $\Hol(\widetilde\fa,\La)$ of the leaf
$\La\backslash\sing(\widetilde\fa)$ of $\widetilde\fa$ is solvable.
By the non-resonant hypothesis  this group contains some element of
the form $f(z) = e^{2\pi i\la}\,z +\dots$ with $\la \in {\bf
C}\backslash\bq$. By Lemma~\ref{Lemma:linearization} this map $f$ is
analytically linearizable and in this same coordinate  the group
$\Hol(\widetilde\fa,\La)$ is either abelian or  it is analytically
conjugate to a subgroup of the group
$$
\bh_k = \big\{\vr(z) = \frac{az}{\sqrt[k]{1+bz^k}}\,, a \ne 0\big\}
\quad\text{for some}\quad k \in \bn.
$$

\noindent In the abelian case the foliation $(\widetilde\fa$ and
therefore) $\fa$ is given by a  closed meromorphic one-form with
simple poles (see \cite{C-LN-S}), defined in a neighborhood of the
origin and therefore it writes as a logarithmic foliation say
$$
\fa:\, \omega = \sum_{j=1}^t \la_j\,\frac{df_j}{f_j}\,, \,\,\, \la_j
\in {\bf C}\backslash\{0\}, \,\,\, f_j \in \O_2\,.
$$
Suppose that $\Hol(\widetilde\fa,\La)$ is solvable and nonabelian.
Then by the main result of \cite{Camacho-Scardua} the foliation
$\fa$ is the pull-back by some germ of a holomorphic map
$\sigma\colon {\bf C}^2,0 \to {\bf C}^2,0$ of a germ of meromorphic
Riccati (Bernoulli type) foliation say
$$
\eR:\, \al(x)dy - \big(y^2\be_0(x) + y\be_1(x))dx=0.
$$
This proves Theorem~\ref{Theorem:germs}.\end{proof}

In the same line of reasoning we can prove:

\begin{Theorem}
\label{Theorem:globalaffine}  Let $\fa$ be a foliation on ${\bf C
P}^2$.  Suppose that:
\begin{itemize}

\item[\rm(i)] $\fa$ is transversely projective
in the complement of an algebraic invariant curve $\Lambda \subset
{\bf C P}^2$.

\item[\rm(ii)] Each singularity $p\in \sing(\fa) \cap \Lambda$ is a
non-degenerate {\rm(}non-dicritical{\rm)} irreducible singularity.

\end{itemize}

\noindent Then $\fa$ is given by a logarithmic one-form in a
neighborhood of the origin or it is a meromorphic pull-back of a
germ of a Bernoulli type foliation $\eR\colon \al(x)dy -
(y^2\be_0(x) + y\be_1(x))dx = 0$ where $\al$, $\be_0$, $\be_1$,
$\be_2$ are meromorphic in some neighborhood of the origin.
\end{Theorem}

\begin{proof}
By Theorem~\ref{Theorem:solvable} the holonomy group of the leaf
$L=\Lambda \setminus (\sing(\fa)\cap \Lambda)$ is solvable. We claim
that some singularity in $\Lambda$ is non-resonant. Indeed, by the
Index theorem (\cite{Camacho-Sad}) the sum of all indexes of
singularities in $\Lambda$ is equal to a (natural) positive number,
the square of the degree of $\Lambda$. This implies that not all
indexes are rational negative. Therefore, the holonomy group of (the
leaf contained in) $\Lambda$ contains some non-resonant germ and we
may proceed as in the proof of Theorem~\ref{Theorem:germs} and apply
the main results in \cite{Camacho-Scardua} to conclude that $\fa$ is
either given a by logarithmic one-form or by a rational pull-back of
a Bernoulli type foliation.
\end{proof}

\begin{Remark}{\rm
 (1) Theorems~\ref{Theorem:germs} and \ref{Theorem:globalaffine} above show
 that in order to capture the generic foliations in the class of
 Riccati foliations it is necessary to allow
 dicritical singularities.

 (2) Theorem~\ref{Theorem:germs} completes
 an example given in \cite{Touzet} of
a germ $\fa$ satisfying (i) and (ii) but which is not a meromorphic
pull-back of a Riccati foliation on an algebraic surface. Indeed,
the construction given in \cite{Touzet} exhibits $\fa$  having as
projective holonomy group $G$, i.e., the holonomy group
$G=\Hol(\widetilde\fa, D)$, where $D$ is the exceptional divisor of
the blow-up, a  non-abelian solvable group conjugate to a subgroup
of $\bh_1 = \big\{z \mapsto \frac{\la z}{1+\mu z}\big\}$. Our result
implies that $\fa$ is a meromorphic pull-back of a germ of
meromorphic Bernoulli foliation in a neighborhood of the origin.

\noindent (3) In \cite{Touzet} it is also given an example of a
foliation $\H$ on a rational surface $Y$ such that  $\H$ is
transversely projective on $Y\backslash X$ for some algebraic curve
$X \subset Y$ and such that $\H$ is {\it not\/} birationally
equivalent to a Riccati foliation on $\ov{\bf C} \times \ov{\bf C}$.
Nevertheless, an analysis of the singular set $\sing(\H)$ of $\H$
shows that our non-dicriticalness and nonresonance hypothesis on the
singularities of the foliation are not satisfied. Indeed, the
singularities are obtained as surface desingularization of quotients
of regular  foliations by finite groups with fixed points.}
\end{Remark}

The next lemma will be useful in capturing the Riccati case.

\begin{Lemma}
\label{Lemma:extendsstructure}

Let $\fa$ be a germ of an irreducible singularity at the origin $0
\in {\bf C}^2$. Assume that:
\begin{enumerate}

\item $\fa$ admits a meromorphic projective triple $(\Omega, \eta,
\xi)$.

\item $\fa$ is transversely projective in the complement of its
local separatrices $\sep(\fa,0)$ which is assumed to be given by
$\{x=0\}\cup\{y=0\}$.

\item The foliation $\fa^\perp=\fa_\xi$ is transverse to the   axis $\{y=0\}$.

\item The foliation $\fa^\perp=\fa_\xi$ is transversely projective
in a neighborhood of the origin minus the axis $\{x=0\}$.

\end{enumerate}

Then $\fa^\perp$ is transversely projective in a neighborhood of the
origin. Indeed, the projective transverse structure of $\fa_\xi$
outside the separatrix $\{y=0\}$ extends as a projective transverse
structure for  $\fa_\xi$ to a neighborhood of the singularity.

\end{Lemma}

\begin{proof}

We have several possible cases regarding the singularity of the
foliation $\fa$:

\noindent{\bf Case 1. $\fa$ is non-degenerate non-resonant}: In this
case, because it is transversely projective in the complement of its
set of local separatrices, by \cite{Touzet} (cf.
Theorem~\ref{Theorem:Touzet})  the singularity is analytically
linearizable. Let us therefore write $\Omega=g(xdy - \lambda y dx)$
in suitable local coordinates, for some meromorphic function $g$ and
$\lambda \in \bf C \setminus \mathbb Q$. By the relations in
Proposition~\ref{Proposition:modificationforms} we may assume that
\begin{equation}
\Omega=\frac{dy}{y} - \lambda \frac{dx}{x}, \, \, \,  \eta=h \Omega
= h \left(\frac{dy}{y} - \lambda \frac{dx}{x}\right), \, \, \,
\xi=-dh - \frac{1}{2} h^2 \left(\frac{dy}{y} - \lambda
\frac{dx}{x}\right) + \ell \left(\frac{dy}{y} - \lambda
\frac{dx}{x}\right)
\end{equation}
where
\[
d\left(\sqrt{\ell} (\frac{dy}{y} - \lambda \frac{dx}{x})\right)=0
\]
Since $\lambda\notin\mathbb Q$ the form $\left(\frac{dy}{y} -
\lambda \frac{dx}{x}\right)$ admits no meromorphic first integral
and therefore we must have $\ell = const=c\in \bf C$. Hence we have
\[
\xi=-\frac{2}{h^2 - 2c}dh -  \left(\frac{dy}{y} - \lambda
\frac{dx}{x}\right)
\]
which is a closed meromorphic one-form. Suppose that $c \ne 0$. In
this case the one-form above has simple poles. On the other hand,
because $\fa_\xi$ is transverse to the axis $\{y=0\}$ we conclude
that this is not contained in the polar set of $\xi$ and therefore
we can integrate $\xi$ as
\[
\xi= \alpha \frac{dx}{x}+ df(x,y)
\]
for some holomorphic function $f(x,y)$. Thus we can write
\[
\xi=\alpha \frac{d(x e^{ \frac{1}{\alpha} f})}{x e^{
\frac{1}{\alpha} f}}
\]
This shows that $\xi$ is  a regular (nonsingular) foliation in a
neighborhood of the singularity, {\it a fortiori}, it is
transversely projective in this neighborhood.

Suppose now that $c=0$. In this case we have
\[
\xi=-\frac{2}{h^2}dh -  \left(\frac{dy}{y} - \lambda
\frac{dx}{x}\right)
\]
The axis $\{y=0\}$ is contained in the polar set of
$\xi$. Because $\xi$ is closed, its polar set is invariant by the
foliation $\fa_\xi$. Since by hypothesis $\{y=0\}$ is not
$\fa_\xi$-invariant we may exclude the case $c=0$.

\noindent{\bf Case 2. $\fa$ is non-degenerate resonant}: In this
case the foliation is a meromorphic pull-back of a Riccati foliation
(\cite{Touzet}). We have two possibilities.

\noindent{\bf (i)} If the foliation is transversely affine outside
of its set of local separatrices then it is either analytically
linearizable as $\Omega=g(xdy - \frac{n}{m}ydx)$ for some $n,m \in
\mathbb N$, or (cf. \cite{Berthier-Touzet})  it is analytically
conjugate to the germ of foliation given by
 \[
 \om_{k,l} = k \,
x\,dy + l \, y(1+\frac{\surd\ov{-1}}{2\pi} x^{l}y^k)dx=0
\]
Thus $\fa$ is given by the closed one-form
\[
\Omega_{k,l}:=\frac{1}{x^{l+1}y^{k+1}}\omega_{k,l}
\]

Proceeding as in the previous case we may assume that
\[
\Omega=\Omega_{k,l}, \,  \eta=h \Omega = h \Omega_{k,l}, \,  \, \,
\xi=-dh - \left(\frac{1}{2} h^2 + \ell\right)\Omega_{k,l}
\]
for some meromorphic function $\ell$ satisfying
$d(\sqrt{\ell}\Omega_{k,l})=0$.

Since $\fa$ is not analytically linearizable, it admits no
meromorphic first integral. Therefore, because $\Omega_{k,l}$ is
closed, we must have $\ell=const.$ Thus we may assume that
\[
\xi=-\frac{2}{h^2 - 2c}dh -  \Omega_{k,l}
\]
and therefore $\fa_\xi$ is given by a closed meromorphic one-form.
At this point the argumentation follows as in the preceding
linearizable case. Again we conclude that $\fa_\xi$ admits a
holomorphic first integral and is transversely projective in a
neighborhood of the singularity.

\footnote{Notice that $\Om_{k,l} = g_{k,l}\,dy_{k,l}$
 where $y_{k,l} = \frac{x ^l y ^{k}}{\frac{\surd\ov{-1}k}
 {2\pi}l x ^{l} y^k   \log x-1}$
and $g_{k,l}\,= \, - \frac {(\frac{\surd\ov{-1}l}{2\pi} x^{l} y^{k}
\log x -1)^2}{x^{l -1}  y^{k-1} }$.}

 \noindent{\bf (ii)} Assume now that the foliation is  not
transversely affine outside of the set of local separatrices. Then
it is a meromorphic pull-back of a Riccati foliation and writes as
\[
\Omega= g (dh - (\frac{1}{2}h^2 - R(f))df
\]
for some holomorphic function $f(x,y)$, some meromorphic function
$R(z)$ and some meromorphic function $h(x,y)$. Because $\fa$ is not
transversely affine in the complement of its set of local
separatrices, we have $R(f)\ne 0$. Moreover, if $\ell$ is a
meromorphic function such that $d(\sqrt{\ell} \Omega)=0$ then $\ell$
is constant. This shows that we may assume that  $\xi=df$. Hence
$\fa_\xi$ is transversely projective in a neighborhood of the
singularity.

\vglue.1in

\noindent{\bf Case 3. $\fa$ is a saddle-node}: If  $\fa$  admits no
affine transverse structure outside of the set of local separatrices
then we proceed as in the second case (ii) above. Assume now that
$\fa$ is transversely affine in the complement of the set of local
separatrices. Then the foliation is  analytically conjugated to its
formal normal form (\cite{Berthier-Touzet}), that is, it can be
written in suitable local coordinates as
\[
\Omega=g(x(1 + \lambda y^k)dy - y^{k+1}dx)
\]
for some meromorphic function $g$, some $\lambda \in \bf C$ and $k \in \mathbb N$.

Then we may assume that
\[
\Omega=\Omega_{k, \lambda}:=\frac{1}{1 + \lambda y ^k}{y^{k+1}}dy  -
\frac{dx}{x}, \, \, \eta= h \Omega_{k , \lambda}= h \frac{1}{1 +
\lambda y ^k}{y^{k+1}}dy - \frac{dx}{x},
\]
and
\[
 \xi=-dh - \left(\frac{1}{2} h^2 +
\ell\right)\Omega_{k,\lambda}
\]
for some meromorphic function $\ell$ satisfying
$d(\sqrt{\ell}\Omega_{k,\lambda})=0$.

As it is well-known the germ $\fa$ admits no meromorphic first integral, therefore
because $\Omega_{k, \lambda}$ is closed, we must have $\ell=const.=c$ and then
\[
\xi=-\frac{2}{h^2 - 2c}dh -  \Omega_{k,\lambda}.
\]

This implies again that $\fa_\xi$ is given by a closed meromorphic
one-form and as above the foliation $\fa_\xi$ admits a holomorphic
first integral in a neighborhood of the singularity. In particular,
$\fa_\xi$ is transversely projective in a neighborhood of the
singularity.
\end{proof}

\medskip

Given a foliation $\fa$ on ${\bf C P}^2$, by an {\it algebraic leaf}
of $\fa$ we mean a leaf $L$ of the foliation which is contained in
an algebraic curve in ${\bf C P}^2$. Thanks to the Identity
Principle and to Remmert-Stein extension theorem, a leaf $L$ of
$\fa$ is algebraic if and only if it accumulates only at singular
points of $\fa$. In this case the algebraic curve consists of the
leaf and such accumulation points. The following remark will be
useful:

\begin{Lemma}
\label{Lemma:commonleaf} Let $\fa$ and $\fa_1$ be distinct
foliations on ${\bf C P}^2$. If a leaf $L$ of $\fa$ is also a leaf
of $\fa_1$ then this leaf is algebraic.
\end{Lemma}

\begin{proof}

 Suppose $(x(z), y(z)),\,\,z\in V\subset \bf C$ is a common
solution of the foliations $\fa$ and ${\fa}_1$ on ${\bf C P}^2$ say:
${\fa}$  is given by $\frac{dy}{dx} = \frac{P(x,y)}{Q(x,y)}$ and $
{\fa}_1$ by $\frac{dy}{dx} = \frac{P_1(x,y)}{Q_1(x,y)}$ where $P,Q$
and $P_1,Q_1$ are relatively prime  polynomials. Then we have
\[
 \frac{P(x(z),y(z)}{Q(x(z),y(z))} =
\frac{dy/dz}{dx/dz} = \frac{P_1(x(z),y(z))}{Q_1(x(z),y(z))}
\]
so that $(PQ_1 - P_1 Q)(x(z),y(z)) = 0$. By  hypothesis $PQ_1 - P_1
Q\not\equiv 0$ so that $L$ satisfies the non-trivial  algebraic
equation $PQ_1 - P_1 Q=0$. It follows that $L$ is algebraic.
\end{proof}

\begin{Theorem}
\label{Theorem:projectivesurfacepossibilities}
 Let $\fa$ be a foliation on a projective surface $M$  with a projective
 transverse structure  outside of an  algebraic curve $\Lambda\subset M$.
 Let $(\Omega, \eta, \xi)$ be a rational projective triple defining
 the projective transverse structure
 outside of the curve $\Lambda$. We have the following
 possibilities:
 \begin{enumerate}
 \item $\Lambda$ contains all the non-dicritical separatrices of
 $\fa$ in $\Lambda$.

 \item $\fa^\perp$ coincides with $\fa$.

 \item $\fa$ is transversely affine in $M \setminus
\Lambda^\prime$ for some algebraic invariant curve
 $\Lambda^\prime\subset M$ containing $\Lambda$.

\item The projective transverse structure of $\fa_\xi$ extends to
$M$.

\end{enumerate}
\end{Theorem}

\begin{proof}
We perform the resolution of singularities for $\fa$ in $\Lambda$
and obtain a projective manifold $\tilde M$, a divisor $E=D\cup
\tilde \Lambda$, where $D$ is the exceptional divisor and $\tilde
\Lambda$ is the strict transform of $\Lambda$,  equipped with a
pull-back foliation $\tilde \fa$ with irreducible singularities in
$E$. The foliation $\tilde \fa$ is transversely projective in
$\tilde M \setminus E$. By Lemma~\ref{Lemma:noninvariantextension}
the projective transverse structure of $\tilde \fa$ extends to the
non-invariant part of $D$ so that, for our purposes we may assume
that $D$ is $\tilde \fa$-invariant, though not necessarily
connected. Take a singular point $q\in \tilde \Lambda\cap
\sing(\tilde \fa)$. Suppose that $\tilde \fa$  exhibits some local
separatrix $\Gamma$ through $q$ which is not contained in $E$. If
$\Gamma$ is $\fa_\xi$-invariant then by Lemma~\ref{Lemma:commonleaf}
$\Gamma$ is contained in an algebraic leaf of $\tilde\fa$ not
contained in $E$. This projects onto an algebraic leaf
$\Lambda^\prime$ of $\fa$ not contained in $\Lambda$. The projective
transverse structure of $\fa$ has $\Lambda^\prime$ as a set of fixed
points and therefore $\fa$ is transversely affine in $M \setminus
(\Lambda\cup \Lambda^\prime)$. Assume now that $\Gamma$ is not
$\fa_\xi$-invariant. Then we are in the situation considered in
Lemma~\ref{Lemma:extendsstructure}. By this lemma we conclude that
the projective transverse structure of $\fa_\xi$ extends from a
neighborhood of $q$ minus the local branches of $E$ through $q$ to a
projective transverse structure in a neighborhood of $q$.

\end{proof}

\begin{Remark}[Logarithmic foliations and invariant curves]
{\rm Theorem~A in  \cite{Licanic} gives the following nice
characterization of  logarithmic foliations: {\it Let $\fa$ be a
holomorphic foliation on a compact algebraic surface $X$ and let $S$
be an invariant compact curve by $\fa$. Assume that one of the
following conditions hold: \noindent (i) $Pic(X)$ is isomorphic to
$\mathbb Z$  or \noindent (ii) $Pic(X)$ is torsion free, $H^1(X, C)
= 0$, $S^2 > 0$ and $\sum\limits_{p\in \sing(\fa)- S} BB_p (\fa) >
O$.
Then, if every local separatrix of $\fa$ through any  $p \subset
\Sing(\fa)\cap S$ is a local branch of $S$ and if every singularity
of $\fa$ in $S$ is a generalized curve, then $\fa$ is logarithmic.}

Here, by $B B_p (\fa)$ we mean  the Baum-Bott index associated to
the Chern number $c_1^2$ of the normal sheaf of the foliation
(\cite{[B-B]}). The author also observes that:

\noindent The condition $\sum\limits_{p\in \sing(\fa)- S} BB_p (\fa)
> O$ holds if each singularity of $\fa$ in $X \setminus S$ is
linearly of Morse type (i.e. f is locally given by the holomorphic
1-form $d(xy) + h.o.t.$). This condition also holds when $\fa$ a
has local holomorphic first integral around each point of $X$ which
is not in $S$.

\noindent Let $\fa$ be a holomorphic foliation on ${\bf C P}^2$ of
degree $m$, then $\sum\limits_{p\in \sing(\fa)- S} BB_p (\fa) = (m +
2)^ 2$.

Therefore, the author proves the  following extension of the second
part of theorem~1 in \cite{Cerveau-LinsNeto} to compact complex
surfaces (cf. \cite{Licanic} Proposition 3.1): {\it Let $\fa$ be a
holomorphic foliation on a compact algebraic surface $X$ with
$H^1(X, \bf C) = 0$ and $Pic(X) = \mathbb Z$. Let $S$ be an
invariant compact curve with only nodal type singularities. If
 $\sum\limits_{p\in \sing(\fa)- S} BB_p (\fa)< S^2$,
then $\fa$ is logarithmic.}

By taking a look at the proof given in \cite{Licanic} we  conclude
that the conclusion of
 Theorem~A holds for a foliation $\fa$ on the complex
 projective plane ${\bf C P}^2$ having
 an invariant algebraic curve $S$ such that each
 singularity of $\fa$ in $S$ exhibits no saddle-node in its
 resolution and if $S$ contains each non-dicritical separatrix
 of each singularity of $\fa$ in $S$.

}
\end{Remark}

\begin{Theorem}
\label{Theorem:mainclassification} Let $\fa$ be a foliation on
$\ov{\bf C} \times \ov{\bf C}$ with a projective transverse
structure in the complement of  an algebraic curve $\Lambda \subset
\ov{\bf C} \times \ov{\bf C}$. Assume  that the singularities of
$\fa$ in $\Lambda$ are all non-resonant generalized curves. Then
$\fa$ is a logarithmic foliation or it is a rational pull-back of a
Bernoulli type foliation or it is a rational pull-back of a Riccati
foliation.
\end{Theorem}
\begin{proof}
First we apply Theorem~\ref{Theorem:projectiveplane} in order to
be able to apply Theorem~\ref{Theorem:projectivesurfacepossibilities}.
According to Theorem~\ref{Theorem:projectivesurfacepossibilities} we
have  the following
 possibilities:
 \begin{enumerate}
 \item $\Lambda$ contains all the non-dicritical separatrices of
 $\fa$ in $\Lambda$.

 \item $\fa^\perp$ coincides with $\fa$.

 \item $\fa$ is transversely affine in $M \setminus
\Lambda^\prime$ for some algebraic invariant curve
 $\Lambda^\prime\subset M$ containing $\Lambda$.

\item The projective transverse structure of $\fa_\xi$ extends to
$M$.

\end{enumerate}

In case (1), because the singularities are non-dicritical  we
conclude that the algebraic curve $\Lambda$ contains all the
separatrices through singularities of $\fa$ contained in $\Lambda$.
Since the singularities are assumed to be generalized curves we may
apply \cite{Licanic} and conclude that $\fa$ is  a logarithmic
foliation.

In case (2) we have $\xi \wedge \Omega=0$ so that $d\eta=0$ and
therefore $\fa$ is transversely affine in the complement of the
algebraic invariant curve given by the polar set of $\eta$. In this
case, thanks to the main result in \cite{Camacho-Scardua} and
\cite{Scardua2} the foliation is a logarithmic foliation or it is a
rational pull-back of a Bernoulli foliation.

In case (3) the situation is similar to the one above in case (2)
and the same conclusion holds.

Finally, in case (4) the foliation $\fa_\xi$ is transversely
projective in $\ov{\bf C} \times \ov{\bf C}$. Since this
manifold is simply-connected we conclude that $\fa_\xi$ admits a
rational first integral. By
Proposition~\ref{Proposition:xifirstintegral} $\fa$ is a rational
pull-back of a Riccati foliation. This ends the proof.

\end{proof}

With a very similar (adapted) proof we have the following variant
for foliations in the projective plane:

\begin{Theorem}
\label{Theorem:mainclassificationprojectiveplane} Let $\fa$ be a
foliation on ${\bf C P}^2$ with a projective transverse structure in
the complement of  an algebraic curve $\Lambda \subset {\bf C P}^2$.
Assume that the singularities  of $\fa$ in $\Lambda$ are all
non-resonant extended generalized curves. Then  $\fa$ is a
logarithmic foliation or it is a rational pull-back of a Bernoulli
type foliation or it is a rational pull-back of a Riccati foliation.
\end{Theorem}

\bibliographystyle{amsalpha}

\vglue.1in

\begin{tabular}{ll}
 Bruno Sc\'ardua:  scardua@im.ufrj.br\\
Instituto de  Matem\'atica -  Universidade Federal do Rio de Janeiro\\
Rio de Janeiro - RJ,   Caixa Postal 68530\\
21.945-970 Rio de Janeiro-RJ \\
BRAZIL
\end{tabular}

\end{document}